\documentclass[review,3p,times]{elsarticle}
\usepackage{amsmath,amsthm,verbatim,amssymb,amsfonts,dsfont,mathrsfs,amscd,graphicx}

\usepackage[utf8]{inputenc}
\usepackage{longtable}
\usepackage{pgfplots}
\usepackage{multirow}
\usepackage{makecell}
\usepackage{enumitem}
\usepackage{xfrac}
\usepackage[title]{appendix}
\usepackage[english]{babel}
\usepackage{pdflscape}
\usepackage{hyperref}
\usepackage{tikz}
\usetikzlibrary{positioning}

\usepackage[nodisplayskipstretch]{setspace}

\usepackage[tableposition=top]{caption}  
\usepackage{subcaption}
\usepackage{tikz}
\usetikzlibrary{arrows,shapes,positioning,shadows,trees,calc}

\usepackage[linesnumbered,ruled]{algorithm2e}

\SetCommentSty{mycommfont}

\def\R{\mathbb R}

\def\N{\mathbb N}
\def\NN{\mathcal N}

\def\X{\mathbb X}

\def\U{\mathcal U}
\def\F{\mathbb F}

\def\L{\mathbb L}
\def\LL{\mathcal L}
\def\LLL{\mathscr L}

\DeclareMathOperator*{\argmin}{arg\,min}

\pgfplotsset{compat=1.17}

\journal{Information Fusion}

\begin{document}

\begin{frontmatter}

\title{Variable Landscape Search: A Novel Metaheuristic Paradigm for Unlocking Hidden Dimensions in Global Optimization}

\author[a,c]{Rustam Mussabayev}
\ead{rustam@iict.kz}
\ead{ru.mussabayev@satbayev.university}

\affiliation[a]{organization={AI Research Lab, Satbayev University},
            addressline={Satbaev str. 22}, 
            city={Almaty},
            postcode={050013},
            country={Kazakhstan}}

\affiliation[b]{organization={Department of Mathematics, University of Washington},
            addressline={Padelford Hall C-138}, 
            city={Seattle},
            postcode={98195-4350}, 
            state={WA},
            country={USA}}

\affiliation[c]{organization={Laboratory for Analysis and Modeling of Information Processes, Institute of Information and Computational Technologies},
            addressline={Pushkin str. 125}, 
            city={Almaty},
            postcode={050010}, 
            country={Kazakhstan}}

\author[a,b]{Ravil Mussabayev\corref{cor1}}
\ead{ravmus@uw.edu}
\ead{r.mussabayev@satbayev.university}

\cortext[cor1]{Corresponding author.}

\begin{abstract}
This paper presents the Variable Landscape Search (VLS), a novel metaheuristic designed to globally optimize complex problems by dynamically altering the objective function landscape. Unlike traditional methods that operate within a static search space, VLS introduces an additional level of flexibility and diversity to the global optimization process. It does this by continuously and iteratively varying the objective function landscape through slight modifications to the problem formulation, the input data, or both. The innovation of the VLS metaheuristic stems from its unique capability to seamlessly fuse dynamic adaptations in problem formulation with modifications in input data. This dual-modality approach enables continuous exploration of interconnected and evolving search spaces, significantly enhancing the potential for discovering optimal solutions in complex, multi-faceted optimization scenarios, making it adaptable across various domains. In this paper, one of the theoretical results is obtained in the form of a generalization of the following three alternative metaheuristics, which have been reduced to special cases of VLS: Variable Formulation Search (VFS), Formulation Space Search (FSS), and Variable Search Space (VSS). As a practical application, the paper demonstrates the superior efficiency of a recent big data clustering algorithm through its conceptualization using the VLS paradigm.

\end{abstract}


\end{frontmatter}

\section{Introduction}

\textit{Metaheuristic methods} are higher-level heuristics designed to select and modify simpler heuristics, with the goal of efficiently exploring a solution space and finding near-optimal or optimal solutions~\cite{Sörensen2018791}. They incorporate strategies to balance local search (exploitation of the current area in the solution space) and global search (exploration of the entire solution space), in order to escape from local optima and reach the global optimum~\cite{Cuevas2021249}. Incorporating a mechanism to adaptively change the search space, the proposed Variable Landscape Search (VLS) approach can be classified as a metaheuristic method, given that it embodies these core elements of metaheuristic optimization strategies.

\textit{Search space} is the space of all potential solutions for a problem. Since the goal of a typical optimization problem is to find a global optimum, and any global optimum is necessarily a local optimum, we can restrict our notion of search space to include only all the local optima for the current objective function landscape~\cite{Locatelli2016251}. In fact, any feasible solution can be mapped to the nearest locally optimal solution via the appropriate local search~\cite{Vanneschi202313}.

A local search provides only a greedy steepest descent to the nearest local optimum from an initial solution. Unless a way to globalize and continue the local search beyond suboptimal local extrema is found, the search will inevitably get stuck in poor solutions~\cite{Boyan2000-learning}. Metaheuristic approaches tackle this core challenge by developing innovative rules and meta-strategies that effectively and efficiently navigate the solution space, ensuring a thorough exploration~\cite{Hussain2019-metaher}. However, most known metaheuristics operate within a static solution landscape, overlooking the potential to leverage information from related, neighboring objective function landscapes. Thoroughly perturbing the original search space and combining the resulting configurations can effectively steer the global search towards more globally optimal solutions.

In the context of an optimization problem, there are two primary ways to influence the search space: altering the problem formulation, modifying the input data, or doing both simultaneously.

Let’s consider both of these possibilities in more detail:

\begin{enumerate}
    \item \emph{Altering the problem formulation:} This involves altering the formulation or structure of the problem to another that is either very similar or equivalent to it~\cite{Mladenović2022405}. The formulation of a problem includes its objectives, constraints, and the relationships between different variables. Tweaking any of these aspects can lead to a different set of solutions that are considered viable. For instance, choosing a number of centroids in the MSSC problem can be considered as changing the problem formulation, thereby altering the search space of the MSSC problem. Similarly, switching between the Monge and Kantorovich formulations of the optimal transport problem can become the source of new interrelated search spaces since the solutions of these two formulations generally may not coincide~\cite{Peyre2019-opttrans}.

    \item \emph{Modifying the input data:} The input data are the raw details, observations, or parameters on which the problem is based. Altering these inputs can change the landscape of the objective function, thereby the distribution of optimal solutions as well~\cite{Audet2017247}. For example, in a resource allocation problem, changing the quantity of available resources would directly impact the resulting optimal solutions~\cite{Dehnokhalaji20171279}. The input data points for the minimum sum-of-squares clustering (MSSC) problem heavily influence the resulting centroids~\cite{Gribel2019-hgmeans}. The choice of input marginal probability distributions or a cost function in the optimal transport problem clearly changes the landscape of optimal transportation plans~\cite{Peyre2019-opttrans};
\end{enumerate}

The input data of an optimization problem can be considered as a subset of a larger, overarching data space, which may include all possible data points that could be considered for problems of this type~\cite{Mussabayev2023-parbigmeans}. One way to modify the input data is by taking a new random sample of data points from this vast pool at each subsequent iteration. Although the resulting sub-sample is randomly composed, it still retains certain common properties with the original data pool, thereby approximating its spatial properties with a smaller number of objects. Additionally, it is worthwhile to explore other methods for altering input data. For instance, introducing random noise by adding new random data points or implementing small random shifts in the positions of the original data points are viable techniques. There exist numerous other potential methods, underscoring the versatility of the proposed approach. Conceptually, systematic minor changes to the input data can reshape the landscape of the objective function and alter the distribution of locally optimal solutions within it, thereby influencing the search process ~\cite{Mussabayev2023}.

Similarly, the subsequent problem reformulation is chosen from a set of all possible ways to frame or define the problem in an equivalent or slightly different manner. This set includes every conceivable method to state the objectives, constraints, and relationships for problems of a similar nature. A systematic slight modification of the problem formulation within a problem formulation space can be organized to yield solution spaces with a preserved degree of commonality, enabling iterative deformation and intentional reshaping of the solution space in a controlled manner to escape local minima~\cite{Pardo2013-vfs,Brimberg2023-vns}. The following methods can be used to change one formulation to another: coordinate system change, linearization, convexification methods, continuous/discrete reformulations, and others~\cite{Mladenović2022405}.

In the paradigm of Variable Landscape Search (VLS), the search space can be understood as the set of all locally optimal solutions that satisfy the constraints and objectives as defined by the slightly varied problem formulation, using the slightly varied input data. This space is critical in optimization as it contains every solution that could potentially be the ``best'' or ``optimal'' solution, depending on the problem's objectives~\cite{Hertz2008-vss}. In the context of VLS, a significant emphasis is placed on the utilization of a diverse array of problem reformulations and varied input data. Despite the apparent heterogeneity of these elements, it is crucial to recognize that the resultant search spaces share more similarities than differences. This commonality, predominantly in the form of similarity among the search spaces, underpins the feasibility of conducting a global search across these spaces. The phenomenon of similar resultant search spaces is essential for enabling the global search process in VLS. If these spaces did not maintain a fundamental similarity, the global search process would be untenable.

Thus, the input data and problem formulation together shape the objective function landscape. This gives rise to an adaptive process that can in systematic manner dynamically combine varied input data and problem formulations, or alter them separately, yielding a sequence of search spaces, each containing distinct yet interconnected locally optimal solutions. Using a suitable local search to explore these dynamic search spaces is expected to find a better solution than searching a static space or only changing the problem formulation.

The systematic processes of problem reformulation and input data alteration can be defined by a structure that provides an order relation among various reformulations and data variations. The concept of neighborhood structures is central to VLS, defining how potential solutions are explored in relation to the current solution. By iteratively and systematically exploring these neighborhoods, VLS increases the effectiveness of the search and promotes comprehensive exploration of the solution space.

Let us define and examine some essential concepts that will be extensively used throughout the article.

A \textit{neighborhood solution} is a solution lying within the defined ``neighborhood'' of a given current solution in the solution space. The concept of a ``neighborhood'' typically relates to a defined proximity or closeness to a given solution~\cite{Altner2014339}. The determination of what constitutes a ``neighborhood solution'' depends on the structure of the problem and the optimization algorithm being used. The process of iteratively exploring these ``neighborhood solutions'' and moving towards better solutions forms the core of local search algorithms and many metaheuristic methods~\cite{Lourenço2019129}. The aim is to gradually improve upon the current solution by exploring its immediate vicinity in the solution space.

A \textit{neighborhood structure} refers to a method or function that can form a set of potential solutions that are adjacent or ``near'' to a given solution in the solution space. The concept of ``nearness'' or adjacency is defined based on the problem at hand and the specific algorithm being used~\cite{Altner2014339}. For instance, in a combinatorial optimization problem, two solutions might be considered neighbors if one can be reached from the other by a single elementary operation, such as swapping two elements in a permutation. This is an example of a trivial neighborhood structure. A more advanced neighborhood structure could involve generating new neighborhood solutions by integrating the current solution with those from its history, i.e., from the sequence of previous incumbent solutions. This method leverages the accumulated knowledge of previous solutions to explore more sophisticated solution spaces.

To solve a specific optimization problem, it is necessary to identify the most relevant neighborhood structure of potential solutions. This structure should facilitate an effective search for the optimal solution. Achieving this requires a deep understanding of the problem, consideration of its specific features, and a priori knowledge that enables efficient transitions from one potential solution to another.

Neighborhood structure is a fundamental concept in many local search and metaheuristic algorithms, such as variable neighborhood search (VNS)~\cite{Hansen2017-vns}, simulated annealing~\cite{Siddique2016}, tabu search~\cite{Gendreau201937}, and genetic algorithms~\cite{Kramer2017}. These methods rely on iteratively exploring the neighborhood of the current solution to try to find a better solution. The definition of the neighborhood structure directly impacts the effectiveness of the search process, the diversity of solutions explored, and the ability of the algorithm to escape local optima. For example, in VNS, the order and rule of altering neighborhoods throughout the search process are significant factors to build the efficient heuristic~\cite{Duarte2018341}.

Assume that for some variation of input data and problem formulation, the optimization task became to minimize objective function $f$ over feasible solution space (feasible region) $S$. Then, the corresponding \textit{fixed objective function landscape} is the fixed set
$$
L_f^S = \{ (x, f(x)) \ | \ x \in S \}
$$
of all feasible solutions paired with the corresponding objective function values that remains constant throughout the optimization process. This set coincides with the graph of the objective function when the function's domain is restricted to the feasible region.

From now on, the terms ``landscape'' or ``objective landscape'' will simply refer to a fixed objective function landscape. Also, we will use the simplified notation $L$ to refer to an arbitrary landscape when its underlying objective function and feasible space are not important for the context.

For a given landscape $L_f^S$, its \textit{search space} $S^*$ is defined to be the subset
$$
S^*\left(L_f^S\right) = \left\{ x \in S \ | \ x \ \text{is a local optimum} \right\} \subseteq S
$$
of all the local optima of the objective function contained in the feasible space. This static search space is characterized by its unchanging nature, where every possible locally optimal solution is defined at the outset of the optimization procedure. All solutions in this space are accessible to be tested, evaluated, or used at any point during the problem-solving process, and no new solutions are introduced or existing ones excluded once the optimization process has commenced.

In the context of optimization algorithms, a \textit{variable objective function landscape} refers to a landscape in which both the distribution of objective function values over the feasible region and the feasible region itself can vary throughout the optimization process due to slight variations in problem formulation and input data. As a result, the landscape's search space is changing as well. Unlike a fixed landscape, where all locally optimal solutions are predefined and static, a variable landscape allows for a degree of modification in the set of potential solutions. Additionally, it can enable the introduction of new potential solutions, or the exclusion of existing ones, based on certain conditions or criteria during the optimization iterations.

The manipulation of a landscape could be based on a variety of the strategies listed above. The primary goal of this controllable dynamism is to enhance the search process, promote diversity, prevent premature convergence to local optima, and increase the likelihood of achieving the global optimum. The concept of a variable objective function landscape is at the core of the Variable Landscape Search (VLS) metaheuristic, where the algorithm systematically modulates the landscape by utilizing some special approach (neighborhood structure) at each iteration, aiming to bypass local minima and aspire towards the global optimum.

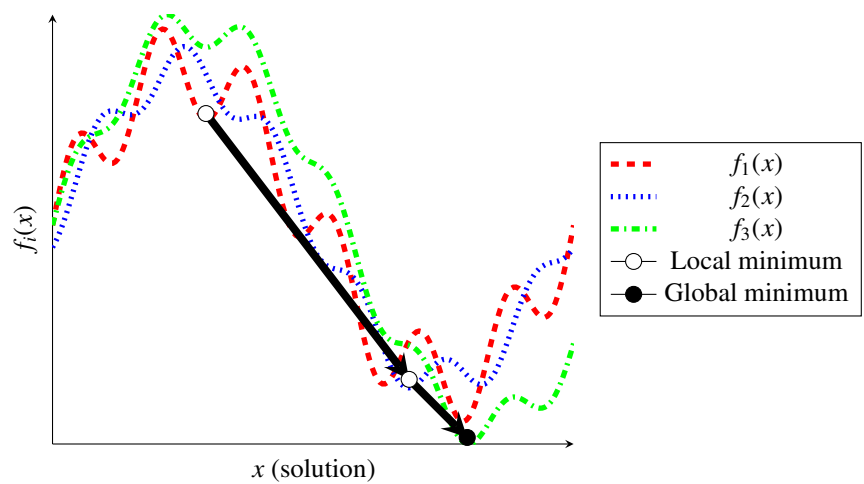
\begin{figure}[htb]
\centering
\begin{tikzpicture}
  \begin{axis}[
    axis lines = left,
    xlabel = $x$ (solution),
    ylabel = {$f_i(x)$},
    legend pos=north west,
    domain=0:4*pi,
    samples=200,
    xtick=\empty,
    ytick=\empty,
    legend style={at={(1.05,0.5)},anchor=west}
    ]
    
    \addplot [red, thick, dashed, line width=2pt]{0.7*sin(0.5*deg(x)) + 0.2*sin(deg(3*x))};
    \addlegendentry{$f_1(x)$}

    \addplot [blue, thick, dotted, line width=2pt]{0.7*sin(0.5*deg(x)) + 0.1*cos(deg(3.0+3*x))};
    \addlegendentry{$f_2(x)$}

    \addplot [green, thick, dashdotted, line width=2pt]{0.9*sin(0.45*deg(x)) + 0.1*sin(deg(3*x))};
    \addlegendentry{$f_3(x)$}
    
    \addlegendimage{mark=*, mark size=3pt, mark options={fill=white, draw=black}}
    \addlegendentry{Local minimum}

    \addlegendimage{mark=*, mark size=3pt, mark options={fill=black}}
    \addlegendentry{Global minimum}
    
    \addplot [only marks, mark=*, mark size=3pt, mark options={fill=white, draw=black}] coordinates {(3.7,0.5)};
    \coordinate (point1) at (axis cs:3.7,0.5);
    
    \addplot [only marks, mark=*, mark size=3pt, mark options={fill=white, draw=black}] coordinates {(8.6,-0.69)};
    \coordinate (point2) at (axis cs:8.6,-0.69);
    
    \addplot [only marks, mark=*, mark size=3pt, mark options={fill=black}] coordinates {(10.0,-0.95)};
    \coordinate (point3) at (axis cs:10.0,-0.95);
    
    \draw [->, >=stealth, line width=3pt] (point1) -- (point2);
    \draw [->, >=stealth, line width=3pt] (point2) -- (point3);

  \end{axis}

\end{tikzpicture}

\caption{Illustration of the process of landscape variation to overcome local minima.}
\label{fig:landscape_varying}
\end{figure}

Figure~\ref{fig:landscape_varying} provides an illustrative example of how manipulating a landscape according to VSL paradigm can potentially advance the global search process. The figure demonstrates the results of three consecutive iterations, each consisting of a minor perturbation of the objective function landscape, followed by a search for the local minimum, applied to an abstract problem. Each arrow signifies an iteration in which a new perturbed landscape $f_i(\cdot)$ is employed, and the algorithm moves from one solution to a better one for the current landscape. This creates a trajectory, allowing the visualization of how the algorithm is navigating the varying landscape, and how it is making progress towards the global optimum. After each iteration, the landscape of the objective function is moderately modified (to a necessary and sufficient extent), which simultaneously extracts the current solution from the current local minimum trap and preserves sufficient commonality among the spatial forms of the resulting landscapes. Through systematic alteration of the objective function landscapes according to the VLS paradigm, an additional ``degree of freedom'' emerges in the optimization process. This allows the search to explore ``additional dimensions'' of the objective function. These dimensions are supplementary but maintain the integrity and commonality of the objective function properties, with the variability/commonality ratio being one of the key parameters.

\section{Related Works}

The concept of manipulating the landscape itself is less common in metaheuristics, but some methods can be interpreted as doing so, either directly or indirectly.

\begin{enumerate}
    \item \emph{Hyper-heuristics~\cite{Drake2020405}:} These are heuristics to choose heuristics, which could be seen as a form of manipulating the search space of the landscape. Instead of operating on the solutions directly, they operate on a space of heuristics that generate or improve solutions. The choice of a heuristic can significantly change the set of potential solutions (local optima) that are explored;

    \item \emph{Constructive metaheuristics}, such as \emph{Ant Colony Optimization (ACO)}~\cite{Dorigo2019311} and \emph{Greedy Randomized Adaptive Search Procedure (GRASP)} ~\cite{Resende2019169}, can be interpreted as manipulating the search space of the landscape. They construct solutions step-by-step, and the choices made in the early steps change the part of the search space that is explored in later steps;

    \item \emph{Cooperative co-evolutionary algorithms (CCEAs)~\cite{Ma2019421}:} In these methods, the solution space is divided into several sub-spaces, and different search processes or different populations evolve in these different sub-spaces. The results are then combined, which can be seen as manipulating the search space;

    \item \emph{Variable Neighborhood Search (VNS)~\cite{Brimberg2023-vns}:} In VNS, the neighborhood structure itself (i.e., the definition of which solutions are ``neighbors'' and can be reached from a given solution) is varied during the search. This can be interpreted as a manipulation of the search space;

    \item \emph{Genetic Algorithms (GA)} with adaptive representation~\cite{Patnaik201745}: Some versions of GAs use adaptive representation, where the encoding of the solutions (i.e., how solutions are represented as chromosomes) is changed during the search based on some criteria. This can be seen as a manipulation of the search space;

    \item \emph{Variable Formulation Search (VFS)}~\cite{Pardo2013-vfs} is a dynamic, adaptive representation of problem formulations rather than a static, predefined one. It refers to the flexibility in the formulation of an optimization problem where the structure of the mathematical representation itself can be altered over the course of the problem-solving process. Unlike a static formulation space, where the mathematical model, decision variables, constraints, and objective function are fixed from the start, a variable formulation space allows for adaptive changes in the problem's formulation. This can involve the introduction of new decision variables, the removal or modification of constraints, or even changes in the objective function, based on certain criteria, algorithmic processes, or evolving insights into the problem during the course of the optimization iterations.

    \item \emph{Formulation Space Search Metaheuristic (FSS)}~\cite{Mladenović2022405} involves a broader strategy where the entire set of possible formulations is considered and structured with a metric or quasi-metric relationship. FSS does not merely adapt a single formulation but explores a structured space of multiple formulations, potentially switching between them based on their relational metrics. This approach extends the search space beyond traditional variable adjustments to include a comparative analysis of different formulations, thus offering a structured and ordered way to navigate through formulation changes.
\end{enumerate}

Data consists of factual elements, observations, or raw information collected through various means, all of which require processing to address specific problems. In the context of optimization, problem formulation generally involves the mathematical representation of a problem, including decision variables, constraints, and objective function(s). These formulations are created independently of any specific data but require data to function effectively. Often, problem formulations are centered around the concepts of functionality, structure, and purpose. For instance, in clustering, the formulation might detail the use of a specific distance metric, the number of clusters to be formed, and any constraints on these clusters~\cite{Correa-Morris20132548}. However, the specific data points to be clustered are inputs to the problem defined by this formulation, not components of the formulation itself.

Problem formulations plays a crucial role in algorithm development by providing a clear framework that specifies the desired outcomes, operational constraints, and conditions relevant to the task. The manner in which a problem is defined and structured profoundly shapes the strategy for designing algorithms, influencing their complexity and the efficacy of the solutions they provide.

The distinction between problem formulation and data is evident, as each serves a unique, separate role. This underscores the importance of addressing them individually when developing novel heuristics. Treating these as separate yet interconnected concepts enhances the clarity and effectiveness of heuristic development, opening new possibilities and dimensions for optimization. This approach allows each conceptual component to be optimized based on its unique requirements and characteristics.

Among all the metaheuristics considered, the ones most similar to our proposed approach are Variable Formulation Search (VFS)~\cite{Pardo2013-vfs} and Formulation Space Search Metaheuristic (FSS)~\cite{Mladenović2022405}. Variable Landscape Search represents a broader strategy in optimization, extending the scope beyond what is offered by FSS and VFS. While the latter focus on exploring different objective landscapes by varying the problem's formulation, Variable Landscape Search extends this concept by not only considering changes in problem formulation but also incorporating modifications in the input data. This dual approach allows for a more comprehensive exploration of potential solutions, making Variable Landscape Search a generalization of VFS and FSS. This approach acknowledges that altering the input data is as crucial as changing the formulation for discovering diverse objective landscapes and their corresponding search spaces. We will explore the details of how Variable Landscape Search generalizes Variable Formulation Search and Formulation Space Search Metaheuristic in Section~\ref{sec:examples}. This section will provide concrete examples of applying Variable Landscape Search to construct optimization heuristics.

\section{Landscape Meta-Space}

The concept of \textit{landscape meta-space} encompasses the aggregation of all possible objective function landscapes relevant to a specific problem or set of problems. Within this framework, each individual objective landscape constitutes a distinct search space --- the set of potential solutions in the form of local optima for that landscape. These potential solutions are defined by a unique problem formulation, representation, constraint set, and the specific dataset being processed.
 
The concept of a landscape meta-space is especially relevant in complex optimization scenarios, where the problem at hand may be addressed from different perspectives or under different assumptions, each leading to a different landscape of the objective function. As such, the landscape meta-space represents the entire universe of potential solutions across all these diverse contexts.

Thus, the landscape meta-space provides a framework for broader exploration and navigation across multiple landscapes, facilitated by a strategically defined neighborhood structure. This structure allows for the shifting or manipulation of the search space itself, rather than merely exploring within a fixed search space. By facilitating transitions between adjacent and potentially more favorable search spaces, it significantly enhances the potential for discovering superior or more diverse solutions.

The landscape meta-space concept takes the exploration process to a higher level of abstraction, allowing search algorithms to consider not only the local optima within a given objective landscape but also the structure and characteristics of the landscape itself. This opens up new possibilities for innovative search strategies, such as the Variable Landscape Search (VLS) method.

\section{Landscape Neighborhood}

A \textit{neighborhood} of a landscape, in essence, encapsulates a subset of the landscape meta-space. This subset can be defined by a `radius of similarity' or `proximity parameter' around the current objective landscape. This `neighborhood' represents a region within the landscape meta-space where each individual landscape shares a certain degree of similarity, correlation, or alignment with the current landscape. These connections may be drawn based on a range of factors, such as solution quality, complexity of problem formulation, or even aspects of data characteristics.

This concept provides a local scope within the broader landscape meta-space, enabling a more focused and computationally efficient search through specially predefined neighborhood structures. The inherent dynamism of a landscape neighborhood is advantageous as it facilitates the exploration of landscapes that are `close' to the current one, potentially leading to improvements in solution quality while maintaining computational feasibility. Through iterative redefinition of the `neighborhood', i.e., by adopting alternative neighborhood structures, the algorithm can progressively navigate the landscape meta-space, yielding a trajectory of landscapes that encapsulates the overall optimization process.

\section{Variable Landscape Search (VLS) Metaheuristic}

For a given optimization problem $P$, there exists the corresponding space of all its formulations $\F$, which is called the formulation space. Each formulation $F = (f, C) \in \F$ consists of a mathematically defined objective function $f$, which is to be minimized, as well as a mathematically defined set of constraints $C$. Alternatively, each set of constraints $C$ can be included into the corresponding objective function $f$ as penalizing terms.

Also, we assume that there exists overarching data space $\X$, which serves as the source of input data for an algorithm solving given problem formulations. Data space $\X$ can be either finite or infinite, representing a data stream in the latter case.

We can influence the objective function landscape $L$ in only two primary ways:
\begin{enumerate}
    \item By changing the input data $X \subseteq \X$;
    \item By changing the problem formulation (or model) $F \in \F$, i.e., by reformulating the problem.
\end{enumerate}

Thus, we have the landscape evaluation map
$$
\LLL: (X, (f, C)) \mapsto L_f^S,
$$
where $X \subseteq \X$ is an input dataset, $F = (f, C) \in \F$ is a formulation for the given problem $P$, and $S$ is the feasible region resulting from evaluating constraints $C$ on input data $X$.

In the most comprehensive implementation of the proposed approach, the search can be conducted through simultaneous variation of the two available modalities $F$ and $X$. Alternatively, each modality can be fixed, allowing searches to be carried out using only one of the available modalities. For example, as a special case, the formulation space $\F$ can be limited to a single fixed formulation $F$. Similarly, the overarching data space $\X$ can be restricted to one fixed dataset $X$, which is then fully utilized at each iteration of the proposed approach. 

The general form of optimization problem that we are looking at may be given as follows:
\begin{equation}
\min \left\{ \ f(x) \mid x \in S, \ L_f^S \in \L \ \right\}
\end{equation}
where $\L$, $S$, $x$ and $f$ denote, respectively, the landscape meta-space (the space of landscapes), a feasible solution space, a feasible solution, and a real-valued objective function. We restrict the landscape meta-space to the image of the landscape evaluation map $\LLL$:
$$
\L = \LLL (2^\X \times \F)
$$
Given that duplicates may exist within the landscape meta-space $\L$, we ensure injectivity by assuming that each landscape $L \in \L$ implicitly contains information from its originating input dataset and formulation. More precisely, for every $L \in \L$ there exists the pair $(X, F) = \LLL^{-1}(L)$. 

\subsection*{VLS ingredients}

The following three steps are executed sequentially in each iteration of a VLS-based heuristic until a stopping condition is met, such as a limit on execution time or number of iterations:

\begin{description}
\item[(1)] Landscape shaking procedure by altering the problem formulation, input data, or both;
\item[(2)] Improvement procedure (local search);
\item[(3)] Neighborhood change procedure.
\end{description}

An initial objective landscape $L_f^S \in \L$ and a feasible solution $x \in S$ are required to start the process, and these starting points are typically selected at random from the corresponding supersets. 

Let us outline the steps listed above within the Basic Variable Landscape Search (BVLS) framework, drawing an analogy with the approach defined in \cite{Brimberg2023-vns}. We refer to this proposed framework as `basic' because, although many other variations of the VLS implementation can be introduced, the following represents the simplest and most straightforward approach.

\subsection*{Shaking procedure}

In this work, we distinguish between different types of a neighborhood structure: on the landscape meta-space, which is the product of neighborhood structures on the data and formulation spaces, and on a feasible solution space.

For an abstract space $E$, let $\NN = \{ \NN_{k_{min}}, \ldots, \NN_{k_{max}} \}$ be a set of operators such that each operator $\NN_{k}$, $k_{min} \le k \le k_{max}$, maps a given choice of element $e \in E$ to a predefined neighborhood structure on $E$:
\begin{equation}
\NN_{k}(e) = \{e' \in E \mid \phi(e, e') \le \Phi_{k} \}, \ k = k_{min}, \ldots, k_{max},
\label{eq:nbhd_struct}
\end{equation}
where $\phi(\cdot, \cdot)$ is some distance function defined on $E$, and $\Phi = \{ \Phi_1, \ldots, \Phi_k \}$ are some positive numbers (integer for combinatorial optimization problems).

Note that the order of operators in $\NN$ also defines the order of examining the various neighborhood structures of a given choice of element $e$, as specified in~\eqref{eq:nbhd_struct}. Furthermore, the neighborhoods are purposely arranged in increasing distance from the incumbent element $e$; that is,
\begin{equation} \label{eq:phi_monotonicity}
\Phi_{k_{min}} < \Phi_{k_{min} + 1} < \ldots < \Phi_{k_{max}}
\end{equation}
Alternatively, each operator $\NN_k$ may be constant with respect to element $e$, yielding a neighborhood structure that is independent of the choice of $e$. Also, the neighborhood structures may well be defined using more than one distance function or without any distance functions at all (e.g., by assigning fixed neighborhoods or according to some other rules).

Now, let us assume that a set of operators $\NN^1$ is defined on the data space $\X$, while another set of operators $\NN^2$ is defined on the formulation space $\F$. Their parameters can be conveniently accessed from the following matrix:
$$
K = \{K_j^i\}_{i \in \{1,2\}, j \in \{min, max\}}
$$

Then, the neighborhood structure on the landscape meta-space $\L$ can be defined as the Cartesian product:
$$
\LL = \NN^1 \times \NN^2
$$

A simple shaking procedure consists in selecting a random element either from $\NN_k^1(X)$ or $\NN_k^2(F)$ depending on the current shaking phase $i$ (see Algorithm~\ref{alg:shaking}). This shaking procedure is presented as an illustrative example; however, other variations or similar procedures could be utilized as alternatives, depending on specific requirements or preferences.

\begin{algorithm}
    \SetAlgoLined
    \SetKwFunction{FMain}{Shake\_landscape}
    \SetKwProg{Fn}{Function}{:}{}
    \Fn{\FMain{$L$, $\LL$, $k$, $i$}}{
        $X, F \gets \LLL^{-1}(L)$\tcp*[l]{obtain data and formulation of incumbent landscape}
        \eIf{$i = 0$}{
            $X' \gets \text{Choose} \ X' \in \NN_k^1(X) \ \text{at random}$\tcp*[l]{shake data}
            $F' \gets F$\;
        }{
            $X' \gets X$\;
            $F' \gets \text{Choose} \ F' \in \NN_k^2(F) \ \text{at random}$\tcp*[l]{shake formulation}
        }
        $L' \gets \LLL(X', F')$\tcp*[l]{evaluate new shaken landscape}
        \Return{$L'$}\;
    }
    \caption{Shaking procedure}
    \label{alg:shaking}
\end{algorithm}

\subsection*{Improvement step}

We assume that a constant predefined neighborhood structure $N_S$ exists on each feasible solution space $S$.

Within an objective function landscape, an improving search typically involves an examination of alternate solutions that are `near' to the current solution $x$. Thus, the improving search uses local information obtained from a single neighborhood $N_S(x)$ consisting of surrounding solutions, where $S$ is the feasible solution space of the current landscape. 

This is precisely what the improvement step does in BVLS. Local search typically utilizes either the `best improvement' strategy, which fully explores $N_S(x)$ and updates $x$ with the optimal solution in $N_S(x)$ if it surpasses the current $x$, or the `first improvement' strategy, which updates $x$ with the first solution in $N_S(x)$ that improves upon $x$. The process concludes when no better solution than the current $x$ is discovered in $N_S(x)$.

The pseudocode for the best improvement local search procedure is presented in Algorithm~\ref{alg:bi_loc_search}.

\begin{algorithm}
    \SetAlgoLined
    \SetKwFunction{FMain}{Best\_improvement\_local\_search}
    \SetKwProg{Fn}{Function}{:}{}
    \Fn{\FMain{$x$, $L_f^S$}}{
        \Repeat{$f(x') \le f(x)$}{
            $x' \gets x$\tcp*[l]{remember old solution}
            $x \gets \argmin_{y \in N_S(x')} f(y)$\tcp*[l]{find best improving solution}
        }
        \Return $x'$\;
    }
    \caption{Local search using best improvement}
    \label{alg:bi_loc_search}
\end{algorithm}

\subsection*{Neighborhood change step}

An iteration of VLS starts with the shake operation, which moves the incumbent landscape $L$ to a perturbed landscape $L'$. Then, the improvement step moves $x$ to a local minimum $x'$; that is,

\begin{equation}
f(x') \le \min \{ f(y) : y \in N_S(x') \}.
\end{equation}

At this point, the neighborhood change step takes over to decide how the search will continue in the next iteration of VLS. The routine commonly used is known as the \emph{sequential neighborhood change step}. If $f'(x') < f'(x)$ (the new solution is better than the incumbent with respect to the objective function of the perturbed landscape), then $x \gets x', \ L \gets L'$  (the new solution with its landscape becomes the incumbent), and $k \gets k_{min}$ (the shake parameter is reset to its initial value); otherwise $f'(x') \ge f'(x)$ holds, $k \gets k + 1$ (increment to the next larger neighborhood or return to the initial one, if $k$ exceeds $k_{max}$), and $n \gets n + 1$ (increment the counter of unsuccessful iterations). The pseudocode for the sequential neighborhood change step is shown in Algorithm~\ref{alg:seq_nbhd_change}. Other forms of the neighborhood change step can be also used. These include cyclic, pipe, random and skewed forms~\cite{Hansen2017-vns}.

\begin{algorithm}
    \SetAlgoLined
    \SetKwFunction{FMain}{Neighborhood\_change\_sequential}
    \SetKwProg{Fn}{Function}{:}{}
    \Fn{\FMain{$x$, $x'$, $L=L_f^S$, $L'=L_{f'}^{S'}$, $k$, $k_{min}$, $k_{max}$}}{
        \eIf{$f'(x') < f'(x)$}{
                $x \gets x'$\tcp*[l]{move in solution space}
                $L \gets L'$\tcp*[l]{move in landscape space}
                $k \gets k_{min}$\tcp*[l]{reset shaking power}
            }
            {
                $k \gets k + 1$\tcp*[l]{increase shaking power}
                \If{$k > k_{max}$}{
                    $k \gets k_{min}$\;
                }
            }
    }
    \caption{Sequential neighborhood change step}
    \label{alg:seq_nbhd_change}
\end{algorithm}

\section{Basic Variable Landscape Search}

The basic VLS scheme requires some stopping criterion for the main loop. Usually, this can be either a limit on the maximum number of non-improving iterations, or a time limit $T$.

Basic VLS operates by alternating between two distinct phases: one focused on shaking within the data dimension $\NN^1$ ($i = 1$), and the other on shaking within the dimension of problem formulation $\NN^2$ ($i = 2$). Each phase $i \in \{1, 2\}$ sets its own limit on the maximum number of iterations within the phase: $S_1 \in \N$ and $S_2 \in \N$, respectively. For convenience, these two integer numbers are combined into a single vector $S = (S_1, S_2)$.

Also, the minimum and maximum values of the shake parameter $k$ are unique for each phase. The values stored in the matrix $K \in \R^{2 \times 2}$ are used to define the admissible ranges of shaking power $k$ in the data and formulation phases, respectively.

Before local search can be applied in the perturbed landscape $L'$, the incumbent solution $x$ should be translated onto the feasible solution space of $L'$ using operator $T_{L,L'}$. The set of operators $\left\{ T_{L,L'} \right\}_{L,L' \in \L}$ governing these transitions should be defined beforehand.

The pseudocode for Basic Variable Landscape Search (BVLS) is given in Algorithm~\ref{alg:bvls}. For initialization, the following starting values can be used:

\begin{itemize}
    \item $L \gets $ Select a current landscape from $\L$ (optionally at random);

    \item $x \gets $ Select an initial current feasible solution from the feasible region $S$ of landscape $L$ (optionally at random).
\end{itemize}

\begin{algorithm}
\SetAlgoLined
\SetKwFunction{FMain}{VLS}
\SetKwProg{Fn}{Function}{:}{}
\Fn{\FMain{$x, L, \LL, K, S, T$}}{
    $t \gets 0$\tcp*[l]{global iteration counter}
    $i \gets 0$\tcp*[l]{current phase: input data (0) or formulation (1)}
    \While{$t < T$}{
        $k \gets K_{min}^i$\tcp*[l]{reset shaking power}
        $s_i \gets 0$\tcp*[l]{iteration counter for current phase $i$}
        \While{$s_i < S_i$}{
            $L' \gets \text{Shake\_landscape}(L, \LL, k, i)$\;
            $x \gets T_{L,L'}(x)$\tcp*[l]{project to feasible set of new landscape $L'$}
            $x' \gets \text{Local\_search}(x, L')$\;
            $\text{Neighborhood\_change\_sequential}(x, x', L, L', k, K_{min}^i, K_{max}^i)$\;
            $s_i \gets s_i + 1$\;
            $t \gets t + 1$\;
        }
        $i \gets (i + 1) \ \text{mod} \ 2$\tcp*[l]{proceed to the next phase}
    }
}
\Return{$x$}
\caption{Basic Variable Landscape Search (BVLS)}
\label{alg:bvls}
\end{algorithm}

\section{Variations and Enhancements in VLS} \label{sec:variations}

It is worth noting that the basic VLS framework (referenced in Algorithm~\ref{alg:bvls}) may employ an alternative strategy to switch between phases: thoroughly exploiting one phase before moving to the next, or alternatively, using only one of the phases. This can be achieved by changing the condition of the inner loop at line 8 of Algorithm~\ref{alg:bvls} to checking if the maximum number of non-improving iterations has been exceeded within the current phase instead of the total number of iterations. However, the current version of VLS presented in Algorithm~\ref{alg:bvls} aims to strike a balance between phase switching and exploitation by choosing the appropriate bounds $S_i$ on the iteration counter for each phase $i$. This approach can be viewed as a kind of natural shaking arising from phase changing.

Another way to speed up the metaheuristic, increase its effectiveness and achieve a natural balance in exploiting different phases can be to distribute VLS jobs between multiple parallel workers in some proportion. Then, either a competitive, collective, or hybrid scheme can be used for communication between the parallel workers~\cite{Mussabayev2023-parbigmeans}. Determining the effectiveness and efficiency of these kinds of parallel settings constitutes a promising future research direction.

\section{Analysis of the Proposed Approach} \label{sec:analysis}

Incorporating a mechanism to adaptively change the objective landscape, VLS can be classified as a metaheuristic method, given that it embodies all the core elements of metaheuristic optimization strategies. VLS harmonizes local search (intensive search around the current area in the solution space) and global search (extensive search across the entire solution spaces of various landscapes), in order to escape from local optima and reach the global optimum.

Within the domain of optimization, the Variable Formulation Search (VFS) method shares several characteristics with the proposed Variable Landscape Search (VLS) approach. However, it is crucial to underscore that VFS is in fact a special case nested within the more encompassing and versatile VLS method proposed in this article.

VFS accomplishes variability in the solution space by facilitating a broadened exploration of potential solutions across diverse problem formulations. This method focuses on exploring different problem representations, objective functions, and constraint sets to create a dynamic landscape of diverse solutions.

Nonetheless, the VLS framework takes this concept a step further. In addition to variable problem formulations, the VLS approach incorporates the modification of the task's input data (datasets) themselves. This extra dimension of flexibility allows for dynamic alteration of the landscape's search space, thereby introducing a higher degree of variability and resilience into the solution search process.

Consequently, it is accurate to posit that the VLS approach is a generalization of the VFS methodology. VLS encapsulates the benefits of variable problem formulation inherent in VFS and extends them to include variations in the input data. The resultant enhanced navigability through a larger, more diversified solution landscape heightens the potential for identifying superior solutions and delivering robust optimization results. Hence, the proposed VLS methodology offers a substantial leap forward in the science of optimization, exhibiting higher adaptability and versatility compared to its predecessors.

The Variable Landscape Search (VLS) metaheuristic stands out among other optimization techniques due to its unique approach of directly and adaptively manipulating the search space itself during the search process. This is a distinctive aspect, as most traditional metaheuristics operate within a fixed search space and focus on manipulating the search process, such as modifying the trajectory, introducing randomness, or adaptively tuning parameters.

Here are some aspects that make VLS unique:

\begin{enumerate}
    \item \emph{Dynamic / Variable Solution Space:} In most metaheuristics, the search space is predefined and remains static throughout the search process. In contrast, VLS changes the search space iteratively, which introduces a higher degree of dynamism and adaptability;

    \item \emph{Balance between Exploration and Exploitation:} The VLS algorithm's systematic manipulation of the objective function landscape helps to strike a better balance between exploration (searching new areas of the solution space) and exploitation (refining current promising solutions). This is achieved by altering the search space itself, which could help avoid entrapment in local minima and promote exploration of previously inaccessible regions of the solution space;

    \item \emph{Robustness:} The VLS approach can enhance the robustness of the search process, as the changing search space can provide different perspectives and opportunities to escape from sub-optimal regions. This could potentially lead to more consistent and reliable optimization results;

    \item \emph{Adaptability:} The VLS is not tied to a specific problem domain. The concept of varying the objective landscape can be applied to a broad range of optimization problems, making VLS a flexible and adaptable metaheuristic.

    \item \emph{Complexity Management:} By actively managing the search space, VLS might provide a novel approach to deal with high-dimensional and complex optimization problems where traditional methods struggle;

    \item \emph{Opportunity for Hybridization:} VLS provides an opportunity for hybridization with other metaheuristics, potentially leading to more powerful and efficient algorithms.
\end{enumerate}

In summary, the VLS approach provides a fresh perspective on optimization, offering potential benefits in terms of robustness, adaptability, and ability to manage complex problem domains. This novel approach to directly and adaptively manipulate the solution space could open new avenues in the design and application of metaheuristics.

\begin{algorithm}
\SetAlgoLined
\KwResult{Compute the final centroids $C$ and cluster assignments $Y$ for a dataset $\X$ using the Big-means algorithm.}
\textbf{Initialization:}\\
Initialize all $p$ centroids $C$ as degenerate\;
$\hat{f}\leftarrow\infty$\;
Set iteration counter $t = 0$\;
\While{$t < T$}{
    Draw a random sample $X$ of size $s$ from $\X$\;
    \For{each centroid $c$ in $C$}{
        \If{$c$ is the centroid associated with a degenerate cluster}{
            \text{Reinitialize $c$ using K-means++ on $X$}\;
        }
    }
    Compute new centroids $C_{\text{new}}$ using K-means on $X$ with initial centroids $C$\;
    \If{$f(C_{\text{new}}, X) < \hat{f}$}{
        $C \leftarrow C_{\text{new}}$\;
        $\hat{f} \leftarrow f(C_{\text{new}}, X)$\;
    }
    $t \leftarrow t + 1$\;
}
$Y \leftarrow \text{Assign each point in } \X \text{ to nearest centroid in } C$\;
\caption{Big-means Algorithm for Big Data Clustering}
\label{alg:big_means}
\end{algorithm}

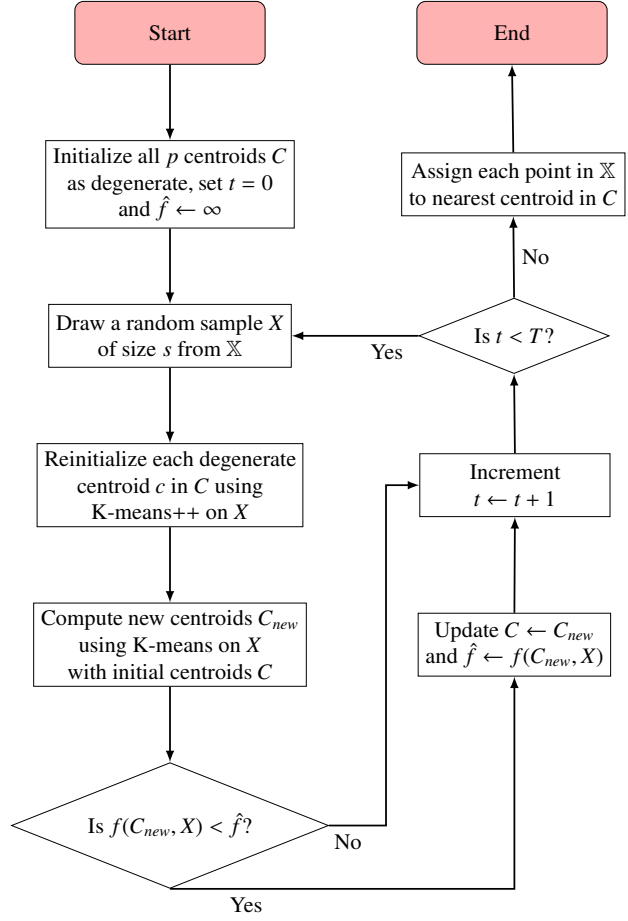
\begin{figure}[htb]
\centering
\resizebox{0.5\linewidth}{!}{
\begin{tikzpicture}[
    node distance=1.2cm and 1.5cm,
    startstop/.style={rectangle, rounded corners, draw, fill=red!30, minimum width=3cm, minimum height=1cm, align=center},
    init/.style={rectangle, draw, minimum width=3cm, minimum height=1cm, align=center},
    operation/.style={rectangle, draw, minimum width=3cm, minimum height=1cm, align=center},
    decision/.style={diamond, draw, minimum width=3cm, minimum height=1.2cm, align=center, inner sep=-2pt, shape aspect=2},
    final/.style={rectangle, draw, minimum width=3cm, minimum height=1cm, align=center},
    arrow/.style={-latex, thick}
]

    \node[startstop] (start) {Start};
    \node[init, below=of start] (init) {Initialize all $p$ centroids $C$ \\ as degenerate, set $t=0$ \\ and $\hat{f} \leftarrow \infty$};
    \node[operation, below=of init] (sample) {Draw a random sample $X$ \\ of size $s$ from $\X$};
    \node[operation, below=of sample] (reinit) {Reinitialize each degenerate \\ centroid $c$ in $C$ using \\ K-means++ on $X$};
    \node[operation, below=of reinit] (new_centroids) {Compute new centroids $C_{new}$ \\ using K-means on $X$ \\ with initial centroids $C$};
    \node[decision, below=of new_centroids, minimum height=2.0cm, minimum width=5.0cm] (decision) {Is $f(C_{new}, X) < \hat{f}$?};
    \node[startstop, right=of start, xshift=0.88cm] (stop) {End};
    
    \node[operation, right=of new_centroids, xshift=0.25cm] (update) {Update $C \leftarrow C_{new}$ \\ and $\hat{f} \leftarrow f(C_{new}, X)$};
    \node[operation, right=of reinit, xshift=0.33cm] (increment) {Increment \\  $t \leftarrow t + 1$};
    \node[decision, right=of sample, xshift=0.5cm] (iter_decision) {Is $t < T$?};
    \node[final, right=of init, xshift=0.18cm] (assign) {Assign each point in $\X$ \\ to nearest centroid in $C$};
    
    \draw[arrow] (start.south) -- (init.north);
    \draw[arrow] (init.south) -- (sample.north);
    \draw[arrow] (sample.south) -- (reinit.north);
    \draw[arrow] (reinit.south) -- (new_centroids.north);
    \draw[arrow] (new_centroids.south) -- (decision.north);
    \draw[arrow] (assign.north) -- (stop.south);
    \draw[arrow] (decision.south) -- ++(1.2cm,0) -| node[pos=0.0, below]{Yes} (update.south);
    \draw[arrow] (decision.east) -- ++(0.9cm,0) |- node[pos=0.0, below, xshift=-0.6cm]{No} (increment.west);
    \draw[arrow] (update.north) -- (increment.south);
    \draw[arrow] (increment.north) -- (iter_decision.south);
    \draw[arrow] (iter_decision.west) -- node[pos=0.0, below, xshift=-0.5cm]{Yes} (sample.east);
    \draw[arrow] (iter_decision.north) -- node[midway, right]{No} (assign.south);
\end{tikzpicture}
}
\caption{Flowchart of the Big-Means algorithm}
\label{fig:big_means_flowchart}
\end{figure}

\section{Examples of Practical and Theoretical Applications} \label{sec:examples}

\subsection{Big Data Clustering Using Variable Landscape Search}

In this section, we will consider a practical example of how we have used the principles underlying the Variable Landscape Search metaheuristic to build an efficient and effective big data clustering algorithm, namely the Big-means algorithm~\cite{Mussabayev2023}.

In the original article~\cite{Mussabayev2023} where the Big-means algorithm was proposed, only its description and a comparative experimental analysis were provided. However, this algorithm was not conceptualized within any existing metaheuristic framework. In fact, attempts to theoretically explain the remarkable efficiency of the Big-means algorithm using existing approaches were unsuccessful. This led us to recognize the uniqueness of the approach used and the possibility of its conceptualization and generalization within a new metaheuristic framework, which is the result presented in this article.

The Big-means algorithm is designed for solving the large-scale Minimum Sum-of-Squares Clustering (MSSC) problem in big data conditions. It uses a heuristic approach and focuses on computational efficiency and solution quality by working with a subset of the data in each iteration, instead of the entire dataset. The pseudocode of Big-means is presented in Algorithm~\ref{alg:big_means}, while its flowchart is depicted in Figure~\ref{fig:big_means_flowchart}.

To start, the algorithm randomly creates a sample $X$ of size $s$ from the given dataset $\X$, where $s$ is much smaller than the total number of feature vectors $m$. The initial configuration of centroids, denoted as $C$, is set using the K-means++ algorithm for the K-means clustering of the first sample. As the algorithm progresses through iterations, the incumbent centroids are updated by the clusterization result of each new sample that improves the clustering solution based on the objective function evaluated on the sample. This approach follows a ``keep the best'' principle, ensuring that the best solution found so far is always prioritized.

Big-means addresses degenerate clusters (also called empty clusters) differently than traditional approaches. When all data points initially associated with a cluster are reassigned to other clusters during the K-means process, Big-means reinitializes those empty clusters using the K-means++ algorithm. This introduces new potential cluster centers and enhances the overall clustering solution by providing more opportunities to minimize the objective function.

The Minimum Sum-of-Squares Clustering (MSSC) formulation aims to partition the data $\X \in \R^{m \times n}$ (possibly very large in both $m$ and $n$) into $p$ clusters, such that the sum of squared Euclidean distances between the data points $x \in \X$ and their associated cluster centers (centroids) is minimized. In mathematical terms, if we denote the centroid of the $j$-th cluster $\X_j$ by $c_j \in C$, then we aim to minimize

\begin{equation}
\min\limits_{C} \ \ f\left(C, \X\right)=\sum\limits_{i=1}^m \min_{j=1,\ldots,p} \| x_i - c_j \|^2
\label{eq:mssc}
\end{equation}
where $p$ is the number of clusters and $\| \cdot \|$ stands for the Euclidean norm.

If we restrict ourselves only to the MSSC formulation~\eqref{eq:mssc}, the space $\L$ of all objective landscapes takes the following form:
$$
\L = \left\{ \ \LLL(X, \text{MSSC}) \ | \ X \subseteq \X \ \right\}
$$
The main idea of the Big-means algorithm and its variations is to follow the stochastic approach, according to which the true $\L$ can be reasonably replaced by its stochastic version:
$$
\widetilde{\L} = \left\{ \ \LLL(X^s, \text{MSSC}) \ | \ X^s \sim \U(\X, l), \ s = 0, \ldots, |\X| \ \right\}
$$
where $\U(\X, s)$ denotes the distribution of $s$-sized samples drawn uniformly at random from $\X$. Thus, the search space of each landscape $L = \LLL(X^s, \text{MSSC}) \in \L$ consists of all locally optimal feasible solutions $C$ for the following MSSC optimization problem defined on sample $X^s$:
\begin{equation}
\min\limits_{C} \ \ f\left(C, X^s\right) = \sum\limits_{x \in X^s} \min_{j=1,\ldots,p} \| x - c_j \|^2
\label{eq:mssc_on_x_s}
\end{equation}

The Big-means paradigm~\cite{Mussabayev2023} further restricts the space $\L$ for an appropriate fixed range of the sample size $s$. The hope is that the resulting objective landscapes would have search spaces reasonably approximating the search space of the original landscape $L^* = \LLL(\X, \text{MSSC})$:
$$
\widetilde{\L}^{[s_{min}, s_{max}]} = \left\{ \LLL(X^s, \text{MSSC}) \ | \ X^s \sim \U(\X, s), \ s = s_{min}, \ldots, s_{max} \right\}
$$
Specifically, the assumption is that by properly and swiftly combining promising local optima of objective landscapes in $\widetilde{\L}^{[s_{min}, s_{max}]}$, one can obtain a close approximation to the global optima of the original objective landscape $L^*$.

Then, the following neighborhood structure can defined on the corresponding data space:
\begin{equation}
\NN_k^1\left(X^s\right) = \left\{ X^{s'} \sim \U(\X, s') \ | \ s_{min} \le s' \le s_{max}, \ |s - s'| \le k \right\}, \ k = K_{min}^1, \ldots, K_{max}^1 \label{eq:bm_data_nbhd_struct}
\end{equation}
The used absolute value distance metric in~\eqref{eq:bm_data_nbhd_struct} clearly satisfies the required monotonicity property~\eqref{eq:phi_monotonicity}, producing a sequence of embedded neighborhoods due to the property that uniform samples of larger sizes subsume uniform samples of smaller sizes.

Moreover, the original Big-means algorithm~\cite{Mussabayev2023} posits that it might be enough to assume $s_{min} = s_{max} = s$ for some appropriate fixed choice of the sample size $s = 0, \ldots, m$ with $s \ll m$:
$$
\widetilde{\L}^s = \left\{ \LLL(X^s, \text{MSSC}) \ | \ X^s \sim \U(\X, s) \right\}
$$

Since algorithms based on the Big-means paradigm consider only one problem formulation (MSSC), they do not employ any neighborhood structure on the formulation space.

Various Big-means-based algorithms fall under the following cases of the general VLS framework (Algorithm~\ref{alg:bvls}):
\begin{enumerate}
    \item \textit{Big-means}~\cite{Mussabayev2023} uses the K-means local search procedure and the sequential neighborhood change (Algorithm~\ref{alg:seq_nbhd_change}) with the following choice of parameters:
    $$
    s_{min} = s_{max} = s \ll m
    $$
    $$
    K_{min}^1 = K_{max}^1 = 0
    $$
    $$
    S_1 = 1
    $$
    $$
    S_2 = 0
    $$
    
    \item \textit{BigOptimaS3}~\cite{Mussabayev2024} maintains several parallel workers, each of which follows the basic VLS scheme, using the K-means local search procedure and the sequential neighborhood change (Algorithm~\ref{alg:seq_nbhd_change}) with the following choice of parameters:
    $$
    1 \le s_{min} \le s_{max} \le m
    $$
    \begin{align*}
    K_{min}^1 &= K_{max}^1 = \infty \ (\text{for the first iteration in a phase}) \\
    K_{min}^1 &= K_{max}^1 = 0 \ (\text{for the remaining iterations in the phase})
    \end{align*}
    $$
    S_1 \in \N
    $$
    $$
    S_2 = 0
    $$
    Also, BigOptimaS3 adds a special procedure to end the search process. It picks the landscape sample size $s_{opt}$ that is most likely to yield improvement in the objective function value based on the collected history of improving samples sizes across iterations. Then, BigOptimaS3 realizes the optimal landscape using $s_{opt}$ and selects the centroids of the parallel worker giving the best result in the optimal landscape;
    
    \item \textit{BigVNSClust} uses the K-means local search procedure and the sequential neighborhood change (Algorithm~\ref{alg:seq_nbhd_change}) with the same choice of VLS parameters as Big-means. Additionally, BigVNSClust adds an extra shaking dimension that acts around the incumbent solution within the feasible solution space of every shaken landscape. This kind of shaking happens simultaneously with the input data shaking. In each iteration, BigVNSClust uses the cyclic neighborhood change scheme to transition from one neighborhood structure to another within the current feasible region, which cyclically increments the shaking power parameter in every iteration irrespective of whether the improvement has taken place or not.
\end{enumerate}
In all these algorithms, the local search within the initial feasible solution space $S$ in Algorithm~\ref{alg:bvls} is performed by starting from the distribution of points defined by K-means++ on $S$, while the incumbent solution $x$ ($C$ in the context of MSSC) is used to initialize the local search procedure in all the subsequent iterations.

The conceptualization and execution of the Big-means algorithm are steeped in fundamental optimization principles. The approach encapsulates and harnesses the power of iterative exploration and manipulation of the search space, overcoming the traditional challenges associated with local minima, and relentlessly pursuing the global optimum.

The perturbation introduced into the clustering results through the `shaking' procedure is a pivotal aspect of the Big-means algorithm. Each iteration yields a new sample, creating variability and diversity in the centroid configurations. Each sample serves as a sparse approximation of the full data cloud in the $n$-dimensional feature space. This stochastic sampling offers a pragmatic balance between exhaustive search (which is often computationally infeasible) and deterministic search (which risks entrapment in local minima). By including a diverse set of sparse approximations, the solution space exploration becomes more robust, adaptive, and capable of reaching the global optimum.

The Big-means algorithm is inherently an adaptive search strategy. Instead of maintaining a fixed search space, it allows the search space to evolve dynamically by sampling different parts of the data cloud in each iteration. This leads to an ``adaptive landscape'', a powerful metaphor where the distribution of locally optimal solutions (the search space) can evolve over the course of optimization, much like species evolving in a changing environment.

This dynamism is beneficial on two fronts: it assists in avoiding the entrapment in local optima, and it promotes a robust exploration of the solution space. If the algorithm is stuck in a local optimum with a given sample, a new random sample might change the landscape such that the local optimum becomes a hill, and better solutions appear as valleys.

The visual representation in Figure~\ref{fig:landscape_varying} can be easily interpreted as Big-means iterations, where $x$ is the incumbent set of centroids $C$, and each new landscape $f_i(C)$ is the MSSC objective function restricted to a new sample $f(C, X_i)$. This figure is a testament to the Big-means algorithm's underlying principle of systematically exploring and manipulating the search space according to VLS approach. The interplay of randomness (through sampling) and determinism (through local optimization) in the algorithm provides a potent strategy to tackle the notorious problem of local minima in clustering algorithms.

The source code for the Big-means algorithm, which includes implementations of various parallelization strategies, is available at \href{https://github.com/R-Mussabayev/bigmeans/}{https://github.com/R-Mussabayev/bigmeans/}.

\subsection{Generalization of Alternative Metaheuristics}

The idea of Variable Formulation Search (VFS) proposed in~\cite{Pardo2013-vfs} retains the same steps as the basic VNS metaheuristic~\cite{Brimberg2023-vns} except for using the special \FuncSty{Accept(x,x',r)} procedure in all of them, where $x$ is an incumbent solution, $x'$ is a candidate solution, and $r$ defines a range of considered formulations $\F = \{ F_1, \ldots, F_r \}$ and their corresponding objective functions $\{ f_1, \ldots, f_r \}$. The \FuncSty{Accept} procedure is listed in Algorithm~\ref{alg:accept_procedure}. This procedure checks if the candidate solution $x'$ leads to improvements in the objective functions by iteratively proceeding to the next formulation in case of a tie in the current one and rejecting $x'$ as long as it causes a decrease in any of them. This approach is effective in tackling the issue of a flat landscape. This implies that, when a single formulation of an optimization problem is employed, a large number of solutions that are close to a given solution tend to have identical values in terms of their objective function. This scenario makes it challenging to identify which nearby solution is the better option for further exploration in the search process. However, our VLS metaheuristic subsumes this VFS idea by simply including the \FuncSty{Accept} procedure into the local search and neighborhood change steps.

\begin{algorithm}
    \SetAlgoLined
    \SetKwFunction{FAccept}{Accept}
    \SetKwProg{Fn}{Procedure}{:}{}
    \Fn{\FAccept{$x$, $x'$, $r$}}{
        \For{$i = 0$ \KwTo $r$}{
            $condition1 \gets f_i(x') < f_i(x)$\;
            $condition2 \gets f_i(x') > f_i(x)$\;
            \uIf{$condition1$}{
                \Return True\;
            }
            \uElseIf{$condition2$}{
                \Return False\;
            }
            \Else{
                continue\tcp*[l]{with the next alternative formulation}
            }
        }
        \Return False\;
    }
    \caption{Accept procedure}
    \label{alg:accept_procedure}
\end{algorithm}

The concept of Formulation Space Search (FSS), as discussed in Mladenovic's study~\cite{Mladenovic2007-fss}, can be considered a specific instance within the broader framework of the Variable Landscape Search (VLS) metaheuristic. This interpretation becomes apparent when you observe that FSS represents what happens when VLS is confined to the single ``formulation'' phase. Essentially, FSS is a special case of VLS where the process is limited to just the formulation aspect.

The idea of Variable Search Space (VSS) proposed in~\cite{Hertz2008-vss} also falls under the VLS metaheuristic. One can check that by restricting to the ``formulation'' phase, modifying the neighborhood change step to the cyclic one, and defining the neighborhood structure on the formulation space as
$$
\NN_k^2 = \{ F_k \}
$$
where $k = K_{min}^2, \ldots, K_{max}^2$ is the shake parameter, and $\F = \{ F_1, \dots, F_r \}$ is a finite set of $r$ available formulations for the underlying optimization problem.

\section{Conclusion}

This paper puts forth the novel Variable Landscape Search (VLS) metaheuristic approach that actively manipulates the search space itself during optimization. Unlike traditional methods that use a fixed, static objective function landscape and its search space, VLS dynamically adjusts them over iterations. This flexibility stems from incorporating changes in the problem formulation as well as modifications of the input data. By expanding the scope of exploration, VLS facilitates the discovery of superior solutions. The algorithm balances focused local search with extensive global search across evolving formulations and input datasets. This helps overcome poor local optima that restrict traditional techniques. The conceptualization of the Big-means clustering algorithm, as well as its various more advanced versions and VFS, under the VLS lens offers insights into how shifting objective landscapes assists in bypassing local minima. The unique strategy and broad applicability of VLS across problem domains position it as a versatile addition to the metaheuristic toolkit for global optimization.

Thus, in this article, we propose a metaheuristic through which a practical result is achieved in the form of Big-means --- a highly effective big data clustering algorithm built on its principles. Additionally, a theoretical result is obtained in the form of a generalization of the following three alternative metaheuristics, which have been reduced to special cases of VLS:

\begin{enumerate}
    \item Variable Formulation Search (VFS)~\cite{Pardo2013-vfs};
    \item Formulation Space Search (FSS)~\cite{Mladenovic2007-fss};
    \item Variable Search Space (VSS)~\cite{Hertz2008-vss}.
\end{enumerate}

An intriguing avenue for future research lies in the enhancement of the VLS metaheuristic by incorporating a shaking mechanism for the incumbent solution within its respective landscape. Such an addition, conceptualized as a distinct phase in the VLS process, aims to further refine and diversify the search strategy. This proposed enhancement has the potential to significantly increase the robustness and efficacy of the VLS metaheuristic in exploring complex search spaces.

\section*{Acknowledgements}

This research was funded by the Science Committee of the Ministry of Science and Higher Education of the Republic of Kazakhstan (grant no. BR21882268).

\bibliography{main}

\end{document}